\documentclass{amsart}
\usepackage{amsmath,amsthm,amssymb, cancel}

\theoremstyle{plain}
\newtheorem{theorem}{Theorem}[section]
\newtheorem{lemma}[theorem]{Lemma}

\theoremstyle{definition}
\newtheorem{definition}[theorem]{Definition}
\newtheorem{example}[theorem]{Example}

\theoremstyle{remark}

\title{A short nonstandard proof of the Spectral  Theorem for unbounded self-adjoint operators}
\author{Takashi Matsunaga}
\address{Department of Medical Informatics, Osaka International Cancer Institute, Osaka City, Japan}
\email{matsunaga-ta@nifty.com}

\begin{document}

\maketitle

\begin{abstract}
Using nonstandard analysis, a very short and elementary proof of the Spectral Theorem for unbounded self-adjoint operators is provided.
\end{abstract}

\renewcommand{\thefootnote}{\fnsymbol{footnote}}
\footnotetext[0]{2020 Mathematics Subject Classification. Primary nonstandard analysis 26E35; Secondary   	spectrum, resolvent 47A10.}

\section{Introduction} 
The Spectral Theorem for unbounded self-adjoint operators (STuB) is one of the most fundamental theorems in functional analysis. Proofs in standard mathematics are lengthy and not straightforward.  The goal of this note is to provide a concise and elementary proof of STuB (Theorem \ref{spec}) by nonstandard analysis. 

Historically, Bernstein offered a nonstandard proof of the Spectral Theorem for bounded self-adjoint operators (STB) in \cite{Bernstein}, while Moore provided alternative nonstandard proof of STB via the nonstandard hull construction in \cite{Moore}. Yamashita and Ozawa \cite{Yamashita} established three equivalent definitions for the nonstandard hull of unbounded self-adjoint operators. Recently, Goldbring presented a nonstandard proof of STuB using the projection-valued Loeb measure  in \cite{Goldbring}, inspired by Raab's \cite{Raab} work.     

Following these devlopments, we provide a very short and elementary nonstandard proof of STuB without using the rather advanced projection-valued Loeb measure which is not covered in textbooks on nonstandard analysis.

\section{Reminders}
We collect useful definitions and lemmas needed for the rest of this note.

\begin{definition}
A spectral family is a family of non-decreasing (orthogonal) projections  $E(\lambda) \ (\lambda \in \mathbb{R})$ on a complex Hilbert space $H$ with $\lim_{\lambda \rightarrow -\infty}E(\lambda)x =0$ and $\lim_{\lambda \rightarrow \infty} E(\lambda)x=x$ for $x \in H$.
\end{definition}

\begin{lemma}\label{int}
Suppose that $E(\lambda)$ is a spectral family on a complex Hilbert space $H$ and
$f(\lambda)$ is continuous on $\mathbb{R}$. Let $D= \{x \in H | \int_{-\infty}^{\infty} f(\lambda)^2 d(E(\lambda)x, x) < \infty \}$ (integrator: $d(E(\lambda)x, x)$). Then for $a<b \in \mathbb{R}$, the Riemann-Stieltjes type integral (integrator: $dE(\lambda)x$)
$$ \int_a^b f(\lambda) dE(\lambda)x \ \ (x \in H), \ \ \ \ \ \ \
    Tx := \int_{-\infty}^{\infty} f(\lambda) dE(\lambda)x \ \ (x \in D)$$
exist. $D$ is dense in $H$. If $f(\lambda)$ is real-valued then $T$ is self-adjoint. 

Let $E(\lambda+0)x = \lim_{\mu \downarrow \lambda} E(\mu)x$ (the right-continuos version of $E(\lambda)$). Then
$$\int_{-\infty}^{\infty} f(\lambda) dE(\lambda+0)x=
    \int_{-\infty}^{\infty} f(\lambda) dE(\lambda)x$$
\end{lemma}
\begin{proof}(Sketch)
Given $\epsilon >0$, let the partition $a = s_0<s_1 < \dots <s_m =b$ be fine so that $|f(s_i)-f(s_{i-1})| < \epsilon/2$. For the first integral, we consider the Riemann sum 
$$Rm(\{s_k\}) := \sum_{k=1}^n f(s_k)(E(s_k)-E(s_{k-1}))x.$$
Let another partition $a = t_0<t_1 < \dots <t_n =b$ with $|f(t_i)-f(t_{i-1})| < \epsilon/2$ and the combined partition $a = u_0<u_1 < \dots <u_p =b \ (p < m+n)$. It is easy to check using the fact that $E(\lambda)$ is a spectral family
$$||Rm(\{s_k\})-(Rm(\{u_k\})||^2 < \epsilon^2||x||^2, \ \ \ \ \
   ||Rm(\{t_k\})-(Rm(\{u_k\})||^2 < \epsilon^2||x||^2.$$
Thus, the first integral exists. For the second integral note that 
$$ ||\int_a^b f(\lambda) dE(\lambda)x||^2 = \int_a^bf(\lambda) d(E(\lambda)x, x)$$
because $E(\lambda)$ is a spectral family. Using this equality, it is shown that the second integral converges. For the third assertion, consider $x_n = (E(n)- E(-n))x$ for $x \in H$. Note that $x_n \rightarrow x\ (n \rightarrow \infty)$ and 
$$ \int_{-\infty}^{\infty} f(\lambda)^2 d(E(\lambda)x_n, x_n) =
\int_{-n}^n f(\lambda)^2 d(E(\lambda)x, x)  < \infty.$$
For the fourth assertion, recall that 
$$ T^* x = \int_{-\infty}^{\infty} \overline{f(\lambda)} dE(\lambda)x $$.
\end{proof}

\begin{lemma}\label{ran}
Suppose that $T$ is a self-adjoint operator on a complex Hilbert space $H$. Then for $z \in \mathbb{C} \smallsetminus \mathbb{R}$
$$ R(T-z) = H \ \ \ (R: \rm{range}).$$ 
\end{lemma}
\begin{proof}
Since $T$ is self-adjoint, $(Tx, y)=(x,Ty)$ for $x, y \in D(T)$ and $T$ is closed. Hence, 
$$||(T-z)x||^2 =||(T-\Re z)x||^2+(\Im z)^2||x||^2 \ge (\Im z)^2||x||^2, $$
so that if $(T-z)x_n \rightarrow y \in H$ then $x_n \rightarrow x_0$ for some $x_0 \in H$ and $(T-z)x_0 =y$. Thus $R(T-z)$ is closed in $H$. 
Let $y \in D(T)$ (dense in $H$) and suppose that for all $x \in D(T)$ $((T-z)x, y) =0$. Since $T$ is self-adjoint, we have $(x, (T-\bar{z})y) =0$ so that $y=0$ by the inequality above. This means the desired results.
\end{proof}

\begin{lemma}\rm{(Operational Calculus)}\label{ope}
Suppose that $E(\lambda)$ is a spectral family on a complex Hilbert space $H$ and let $Tx =  \int_{-\infty}^{\infty} \lambda dE(\lambda)x$. Then for $z \in \mathbb{C} \smallsetminus \mathbb{R}$ and $x \in H$
$$ \int_{-\infty}^{\infty}\frac{1}{\lambda-z} dE(\lambda)x = (T-z)^{-1}x$$
\end{lemma}
\begin{proof}
$$ (T-z)\int_{-\infty}^{\infty}\frac{1}{\lambda-z} dE(\lambda)x = \int_{-\infty}^{\infty}(\mu-z) (\int_{-\infty}^{\infty}\frac{1}{\lambda-z} dE(\lambda)x)dE(\mu)$$
$$ = \int_{-\infty}^{\infty}\int_{-\infty}^{\infty}\frac{\mu-z}{\lambda-z}dE(\lambda)dE(\mu)x = \int_{-\infty}^{\infty}\frac{\nu-z}{\nu-z} dE(\nu)x =\int_{-\infty}^{\infty}dE(\nu)x =x. $$
This leads to the desired result.
\end{proof}

\begin{lemma}\label{del}
For $a<b \in \mathbb{R}$
$$ \lim_{\epsilon \downarrow 0}
\frac{1}{2{\pi}i}\int_a^b(\frac{1}{\lambda-\mu-i\epsilon}-\frac{1}{\lambda-\mu+i\epsilon})d\mu
= \ \ 0 (\lambda<a, \lambda>b ), \ 1(a<\lambda<b).$$
\end{lemma}
\begin{proof}
(Sketch) Recall a representation of the Dirac delta function 
$$\delta(\lambda-\mu) = \lim_{\epsilon \downarrow 0}\frac{1}{2{\pi}i}(\frac{1}{\lambda-\mu-i\epsilon}-\frac{1}{\lambda-\mu+i\epsilon}).$$
\end{proof}

\section{Preliminaries}
We assume familiarity with basic nonstandard analysis. A quick introduction to nonstandard analysis is presented in the Appendix. Consult Davis \cite{Davis} for details.

We extend Moore's \cite{Moore} idea to prove Theorem \ref{spec}. Let $T$ be as in Theorem \ref{spec}, and let $S$ be a *-finite(hyperfinite) -dimensional subspace of $^*D(T)$ such that $D(T) \subseteq S$ (the Concurrence Principle). Note that we assume $H \subset \,^*\!H$ here. Let ${\rm fin}(S) = \{x \in S | \ ^*||x|| \rm{\ is \ finite} \}$, and let $\hat{S} = {\rm fin}(S)/\simeq$ (the quotient space), where $ x \simeq y$ means that $^*||x-y||$ is infinitesimal. Let $\pi:{\rm fin}(S) \rightarrow \hat{S}$ be the canonical quotient mapping. Difine $(\pi(x), \pi(y)) \equiv {\rm st}{(^*(x,y))}$. $\hat{S}$ is called the nonstandard hull of $S$.  By completion, if necessary, we can assume $\hat{S}$ is a Hilbert space. Since $D(T)$ is dense in $H$, $H \subseteq \hat{S}$.  If $A$ is a *-linear operator on $S$ such that $^*||A||$ is finite, one can define the nonstandard hull $\hat{A}$ of $A$ on $\hat{S}$ by setting $\pi(Ax) = \hat{A}(\pi(x))$. 

Let $P_S$ denote the *-projection of $^*H$ onto $S$, and let $T_S$ be the restriction of $P_S{^*T}$ to $S$. It is easy to check that $T_S$ is a *-finite-dimensional *-self-adjoint operator. Hence for $x \in S$, $ T_S  x= {^*\sum} \lambda_n P_n  x$ (*-finite sum)
by transferring the finite-dimensional Spectral Theorem in \cite{Halmos}, where $P_n$ is the *-projection corresponding to the eigenvalue $\lambda_n$. For $\lambda \in {^*\mathbb{R}}$ by setting $F(\lambda) = {^*\sum_{\lambda_n \le \lambda}}P_n,$
we obtain a *-spectral family $F(\lambda)$ on $S$.  For a fixed $x \in S$, we can define the *-Riemann Stieltjes integral, and we obtain $T_S x =  {^*\sum} \lambda_n P_n  x = {^*\int_{-\infty}^\infty} \lambda dF(\lambda) x$. Since $F(\lambda)$ is a *-spectral family, if $x \in R$ then $F(\lambda)x \in R$, where $R = \{x \in S \ | \  ^*||T_Sx|| {\rm \ is \ finite} \}$. It is easy to see that $\hat{F}(\lambda) \  (\lambda \in \mathbb{R})$ is a famliy of non-decreasing projections on $\hat{S}$. Let  $\hat{R} = {\rm the \ closure \ of} \  \pi(R)$ in $\hat{S}$. $\hat{R}$ is a closed linear subspace of $\hat{S}$. Thus, $\hat{F}(\lambda)$ is a spectral family on $\hat{R}$.

From now on, we omit  the upper-left * if there is no ambiguity.

\begin{lemma}{\rm (Representation Lemma)}{\label{rep}}
For $x \in R$
 $$ \pi({^*\int_a^b} \lambda dF(\lambda)x) = \int_a^b \lambda d\hat{F}(\lambda)\pi(x). $$
\end{lemma}

\begin{proof}
Let $a = a_0<a_1 < \dots <a_n =b$ with $max_k (a_k-a_{k-1})^2 < \epsilon$. Noting that $F(\lambda)$ is a *-spectral family, we have

$$ ||\pi(^* \int_a^b \lambda dF(\lambda)x) - \sum_{k=1}^n a_k((\hat{F}(a_k)-\hat{F}(a_{k-1}))\pi(x)||^2 $$
$$ = ||\sum_{k=1}^n \pi(^*\int_{a_{k-1}}^{a_k} \lambda dF(\lambda)x) - \sum_{k=1}^n \pi(a_k((F(a_k)-F(a_{k-1}))x)||^2 
 = ||\sum_{k=1}^n \pi(^*\int_{a_{k-1}}^{a_k} (\lambda -a_k) dF(\lambda)x) ||^2 $$
$$   =\sum_{k=1}^n ||\pi(^*\int_{a_{k-1}}^{a_k} (\lambda -a_k) dF(\lambda)x) ||^2 <\sum_{k=1}^n  \epsilon ||\pi(^*\int_{a_{k-1}}^{a_k} dF(\lambda)x) ||^2  = \epsilon ||\pi((F(b)-F(a))x)||^2 \le \epsilon||x||. $$
By letting $\epsilon  \downarrow 0$, the desired result is obtained.
\end{proof}

\section{A nonstandard proof the Spectral theorem}

\begin{lemma}{\label{inf}}
If $x \in R$, in particular, $x \in D(T)$, for $K_0 = ^*\!\!\mathbb{R} \smallsetminus [-K, K] \  (K:  \rm  positive\ infinite)$
$$ ^*\int_{K_0} \lambda dF(\lambda)x \simeq 0, \ \lim_{\lambda \rightarrow -\infty} \hat{F}(\lambda)\pi(x) =0, \ \lim_{\lambda \rightarrow \infty} \hat{F}(\lambda)\pi(x) = \pi(x). $$
\end{lemma} 
\begin{proof}
Since $F(\lambda)$ is a *-spectral family, for $x \in R$
$$ ||^*Tx||^2 \!\ge\! ||P_S{^*Tx}||^2 \!=\! ||T_Sx||^2 \!=\! {^*\!\!\!\int_{-\infty}^{\infty}} {\lambda}^2 d(F(\lambda)x,x) \ge K|\,^*\!\!\!\int_{K_0} \lambda d(F(\lambda)x,x)| .$$
By assumption, $ ^*\!\!\int_{K_0} \lambda d(F(\lambda)x,x)$ is infinitesimal. Using the *-Polarization Identity (note that $S$ is *-finite dimensional and the integral is actually just a *-finite sum), $(^*\int_{K_0} \lambda dF(\lambda)x, y)$ is also infinitesimal for $x, y \in R$, leading to the first formula. In particular, we have $F(-K)x \simeq 0$ and $F(K)x \simeq x$ so that the second and third formulas easily follow by the definition of $\hat{F}(\lambda)$.
\end{proof}

\begin{lemma}\label{genspec}
For $x \in D(T) \subseteq R$
 $$ Tx = \pi (T_Sx) = \pi({^*\int_{-\infty}^\infty} \lambda dF(\lambda)x) = \int_{-\infty}^{\infty} \lambda d\hat{F}(\lambda)x
= \int_{-\infty}^{\infty} \lambda dE(\lambda)x, $$
where $E(\lambda) = \lim_{\mu \downarrow \lambda}\hat{F}(\mu)$.
\end{lemma}
\begin{proof}
The first equality arises from the fact that $Tx \in H$ and $H \subseteq \hat{S}$. The second equality is the Spectral Theorem for the *-finite dimensional self-adjoint operator $T_S$. For the third equality, applying Lemmas \ref{rep} and \ref{inf}, let $a \downarrow -\infty$ and $b \uparrow \infty$. Lemma \ref{int} gives the last equality.
\end{proof}

\begin{lemma}\label{Hspec}
$E(\lambda)x \in H$ for $x \in H$ and it is unique.
\end{lemma}
\begin{proof}
Since for $\mu \ge \lambda$ $||(\hat{F}(\mu)-\hat{F}(\lambda))x||^2 = (\hat{F}(\mu)x, x)-(\hat{F}(\lambda)x, x) \ge 0$, $\hat{F}(\lambda)x$ has at most countable discontinuities in $\lambda$. Thus, if $\hat{F}(\lambda)x$ is continuous in $\lambda$ at $a<b \in \mathbb{R}$, using Lemma \ref{del} to obtain the first equality and Lemma \ref{ope} to obtain the fourth equality and the (fifth) set membership, we have for $x \in H$
\begin{align*}
\hat{F}(b)x-\hat{F}(a)x
&=\frac{1}{2{\pi}i}\int_{-\infty}^{\infty}\lim_{\epsilon \downarrow 0}\int_a^b(\frac{1}{\lambda-\mu-i\epsilon}-\frac{1}{\lambda-\mu+i\epsilon})d\mu d\hat{F}(\lambda)x \\
&= \lim_{\epsilon \downarrow 0} \frac{1}{2{\pi}i}\int_{-\infty}^{\infty}\int_a^b(\frac{1}{\lambda-\mu-i\epsilon}-\frac{1}{\lambda-\mu+i\epsilon})d\mu d\hat{F}(\lambda)x \\
&= \lim_{\epsilon \downarrow 0} 
\frac{1}{2{\pi}i}\int_a^b\int_{-\infty}^{\infty}(\frac{1}{\lambda-\mu-i\epsilon}-\frac{1}{\lambda-\mu+i\epsilon}) d\hat{F}(\lambda)d\mu \, x\\
&= \lim_{\epsilon \downarrow 0} \frac{1}{2{\pi}i}\int_a^b((T-\mu-i\epsilon)^{-1}-(T-\mu+i\epsilon)^{-1})d\mu \, x \ \ \in H, \\
\end{align*}
Letting $a \downarrow -\infty$, we have $\hat{F}(b)x \in H$ using Lemma \ref{inf}. Thus $E(\lambda)x = \lim_{b \downarrow \lambda}\hat{F}(b)x \in H$. The uniquness is clear from the above expressions.
\end{proof}

Lemmas \ref{inf}, \ref{genspec}, and \ref{Hspec} provide Theorem \ref{spec}.

\begin{theorem}{\rm (Spectral Theorem)}{\label{spec}}
Let $T$ be a (possibly unbounded) self-adjoint operator on a dense subspace $D(T)$ of a complex Hilbert space $H$. Then, there exists a spectral family $E(\lambda)$ on $H$ such that for $x \in D(T)$ 
 $$ Tx = \int_{-\infty}^{\infty} \lambda dE(\lambda)x.$$
Moreover, if $E(\lambda)x$ is right-continuous in $\lambda$, it is unique.
\end{theorem}

\section{Concluding Remarks}
For those familiar with nonstandard analysis, this is a very short and elementary proof of the Spectral Theorem for unbounded self-adjoint operators. 

\subsection*{Disclosure statement and funding}
There are no interests to declare. No funding was received.

\section{Appendix: a quick introduction to nonstandard analysis }
In this appendix, we assume that the reader is fimiliar with basics of set theory and first order logic. 

Nonstandard analysis is a theory founded by Abraham Robinson in the 1960s, motivated largely by the revival of Leibnizian infinitesimals. He constructed a proper extension of $\mathbb{R}$ denoted by $^*\mathbb{R}$, which is logically similar to $\mathbb{R}$ but includes ideal elements such as infinitesimals and infinite numbers. In essence, he  constructed a nonstandard extension $^*U (\ni {^*\mathbb{R}})$ (Theorem \ref{extension}) of a universe (Definition \ref{universe}) $U (\ni \mathbb{R})$, where $^*U$ is logically similar to $U$ but includes ideal elements. 

We need several definitions and lemmas to prove Thoerem \ref{extension}.

\begin{definition}
A (logical) formula is defined  recursively as follows: 
\begin{enumerate}
\item For variables or constants $u$ and $v$, $u=v$ and $u \in v$ are (atomic) formulae,
\item If $\phi$ is a formula, $\lnot\phi$, $\exists x \phi$, and $\forall x \phi$ are formulae, 
\item If $\phi$ and $\psi$ are formulae, $\phi \land \psi$, $\phi \lor \psi$, and $\phi \to \psi \equiv \lnot\phi \lor \psi $ are formlae, 
\item Only those obtained by the above rules are formlae.
\end{enumerate}
A bounded variable $x$ in $\phi$ is a variable that appears in the form of  $\exists x \phi$ or $\forall x \phi$. Other variables in $\phi$ are free variables. A sentence is a formula without free variables.
\end{definition}

\begin{definition}\label{universe}
A universe $U$ is a set satisfying the following conditions:
\begin{enumerate}
\item $u \in v$ and $v \in U$ imply $u \in U$, 
\item $u \in U$ and $v \in U$ imply $\{u, v\} \in U$,
\item $u \in U$ implies $\bigcup u \in U$,
\item $u \in U$ implies $\mathcal{P}(u) \in U$, where $\mathcal{P}(u)$ is the power set of $u$. 
\end{enumerate}
\end{definition}

\begin{example}(Universe containing $\mathbb{R}$)
Let $V_0 = \mathbb{R}$ and define $V_{n+1}$ by $V_{n+1} = \bigcup V_n$ inductively. Set $U_0 = \bigcup V_n$. Define $U_{n+1}$ by $U_{n+1} = U_n \cup \mathcal{P}(U_n)$ inductively. Set $U = \bigcup U_n$. It is straightforward to verify that $U$ is a universe.
\end{example}
Note that $U$ contains various $\mathbb{R}$-related objects such as real numbers themselves, functions on $\mathbb{R}$, and binary relations on $\mathbb{R}$ and so on. 
Usually, we adopt a universe $U$ that has all relevant mathematical objects.

\begin{definition}
A filter basis $\mathcal{F}$ on a set $I$ is a subset of $\mathcal{P}(I)$ satisfiying (1) and (2). A filter $\mathcal{F}$ is a subset satisfying (1), (2) and (3). An ultrafilter $\mathcal{F}$ is a subset satisfying (1), (2), (3) and (4). 
\begin{enumerate}
\item $\phi \not\in \mathcal{F}$ and $I \in \mathcal{F}$, 
\item If $A \in \mathcal{F}$ and $B \in \mathcal{F}$, then $A \cap B \in \mathcal{F}$,
\item If  $A \in \mathcal{F}$ and $A \subseteq B$, then $B \in \mathcal{F}$,
\item If $\mathcal{G}$ is a filter and $\mathcal{F} \subseteq \mathcal{G}$, then $\mathcal{F} = \mathcal{G}$, that is, $\mathcal{F}$ is maximal under $\subseteq$.
\end{enumerate}
\end{definition}

\begin{lemma}\label{ultraexist}
For a filter basis $\mathcal{F}_0$ on $I$, there exists an ultrafilter $\mathcal{F}$ containing $\mathcal{F}_0$.
\end{lemma}
\begin{proof}
Let $\mathcal{F}_1 = \{ X \in \mathcal{P}(I) | X \supseteq A \ {\rm for \ some} \ A \in \mathcal{F}_0 \}$. It is easy to verify that $\mathcal{F}_1$ is a filter. If $\mathcal{F}_1$ is not maximal, there exsits a filter $\mathcal{F}_2 \supsetneq \mathcal{F}_1$. If $\mathcal{F}_2$ is not maximal, there exsits a filter $\mathcal{F}_3 \supsetneq \mathcal{F}_2$ and so on. Finally we have a maximal $\mathcal{F}$. More formally, the existence of $\mathcal{F}$ follows by Zorn's lemma. 
\end{proof}

\begin{lemma}\label{either}
For an ultrafilter $\mathcal{F}$ on $I$ and $A \subseteq I$, either $A \in \mathcal{F}$ or   
$I - A \in \mathcal{F}$ holds.
\end{lemma}
\begin{proof}
Suppose that $A \not\in \mathcal{F}$. Let $\mathcal{G}  = \{X \in \mathcal{P}(I) | X \cup A \in \mathcal{F} \}$. It is a routine to check that $\mathcal{G}$ is a filter containing $\mathcal{F}$. Hence $\mathcal{G} = \mathcal{F}$ by hypothesis. Obviously $I-A \in \mathcal{G}$. This completes the proof. 
\end{proof}

\begin{definition}(Ultrapower)
Let $V$ be an infinite set and $I$ be an infinite index set. Let denote the set of all maps from $I$ to $V$ by $V^I$. Let $<a(i)>, <b(i)> \in V^I$. Define the equivalence relation on $V^I$ by $\{i \in I | a(i) =b(i) \} \in \mathcal{F}$ (this is well-defined if $\mathcal{F}$ is an ultrafilter on $I$). The ultrapower of $V$ over an ultrafilter $\mathcal{F}$ on the index set $I$ is the set of equivalence classes of $V^I$.
\end{definition}

\begin{theorem}\label{extension}(Nonstandard Extension)
For a given universe $U$, there exsits a nonstandard extension $^*U$ and a binary relation $^* \! \! \in$ on it that satisfy the following two conditions:
\begin{enumerate}
\item the Transfer Principle: There exists an injective map $^*: U \rightarrow {^*U}$ such that any sentence $\phi$ in $U$ holds if only if the "corresponding" sentence $^*\phi$ holds in $^*U$. Here the "corresponding" sentence $^*\phi$ is defined by replacing $\in$ with $^*\!\in$ and $c$'s (constants) with $^*c$'s in $\phi$,
\item the Concurrence Principle: For a formula $\phi(a,b)$ "concurrent" with respect to $A (\in U)$, then there exists $b \in {^*U}$ such that $^*\phi(\,^*a, b)$ holds for all $a \in A$. Here "conccurent" with respect to $A$ means that for any finite collection $a_i \in A (\in U)$, there exists $b \in U$ such that $\phi(a_i, b)$ holds for all $i$.
\end{enumerate}
\end{theorem}
\begin{proof}
The construction of $^*U$ proceeds as follows. Consider the index set $I$ of all finite subsets of $U$. For each $i \in I$, let $\mu(i) = \{ j \in I \mid i \subseteq j \}$. Since $\mu(i_1) \cap \mu(i_2) \cap \cdots \cap \mu(i_n) = \mu(i_1 \cup i_2 \cup \cdots \cup i_n)$, this family is a filter basis on $I$ so that  there exists an ultrafilter $\mathcal{F}$ on $I$ containing all the $\mu(i)$'s by Lemma \ref{ultraexist}.  Denote the ultrapower of $U$ over $\mathcal{F}$ by $^*U$ and the equivalence class by $[a(i)], [b(i)]$. The binary relation $[a(i)] \, {^*\!\!\in}\ [b(i)]$ is also well-defined by $\{ i \in I \mid a(i) \in b(i) \} \in \mathcal{F}$. The map * is defined by $a \mapsto [a(i) = a]$. This map is injective by definition.

The Transefer Principle is the special case of the next lemma. For the Concurrence Principle, by assumption, there exists $b(i)$ such that $\phi(a, b(i))$ holds for all $a \in i \cap A$, we obtain the desired result again using the next lemma.
\end{proof}

\begin{lemma}(\cancel{L}os's Theorem)
Under the same conditions as in Theorem \ref{extension}, for any formula $\phi(x,y,\dots,z)$, $^*\phi([a(i)], [b(i)], \cdots, [c(i)])$ holds if only if
 
$\{i \in I \ |\ \phi(a(i), b(i), \cdots ,c(i)) {\rm \ holds} \} \in \mathcal{F}$, where $x, y, \cdots, z$ are free variables.
\end{lemma}
\begin{proof}
We apply mathematical induction on the number of logical symbols in $\phi(x,y, \dots, z)$. First, if there is no logical symbols in $\phi(x,y,\cdots, z)$, $\phi(x,y,\cdots, z)$ is $x\in y$ or $x=y$ so that the conclusion follows by definiton. Next, using De-Morgan's law, it suffice to prove the cases  of $\exists, \land$ and $\lnot$. Consider the case of $\exists$. Clearly, $\exists z \phi([a(i)],[b(i)], \cdots, z)$ holds if only if for some ${c(i)}$   $\phi([a(i)],[b(i)], \cdots, [c(i)])$ holds. By the hypothesis of induction, this occurs if only if $\{i \in I \ |\ \phi(a(i), b(i),$

$\cdots, c(i)) {\rm \ holds} \} \in \mathcal{F}$. The other cases are similar and left to the reader. The key to the proofs is Lemma \ref{either}.  
\end{proof}

\begin{definition}(Standard, Internal, External)
An entity $u \in U$ and the correponding $^*u \in {^*U}$ are called standard. An entity $v \in {^*U}$ such that  $v \,^*\!\in {^*u}$ for some standard $u$ is called internal. Otherwise $v$ is called external.
\end{definition} 

\begin{example}(Transfer Principle I: Embedding, Exstension, *-Omission)
\begin{enumerate}
\item If $A \in U$ and $a \in A$, then $a \in U$ by the transitivity of $U$. Thus, $^*a$ is defined and $^*a  {\,^*\!\in} {\,^*\!A}$. Since the map $^*: U \mapsto \,^*U$ is injective, $A$ is embedded into $^*A$. We often assume $A \subseteq {^*A}$, if there is no confusion,
\item For $a, b \in \mathbb{R}$ $a<b$ if only if ${^*a} {\ ^*\!\!<} {\ ^*b}$. If we assume $\mathbb{R} \subseteq {^*\mathbb{R}}$ as in (1), $^*\!\!<$ is regarded as an extension of $<$ so that $^*$ is often omitted,
\item For $A, B \in U$ and $f: A \mapsto B$, $b=f(a)$ if and only if $^*b= {^*f(^*a)}$. If we assume $A \subseteq {^*A}$ and $B \subseteq {^*B}$ as in (1), $^*f$ is again an extension of $f$ so that $^*$ is often omitted.
\item For a function $f$ on $\mathbb{R}$,
$$\forall x \in \mathbb{R} \, \forall y \in \mathbb{R} \ f(x+y) = f(x) + f(y)$$
if only if 
$$\forall x \,^*\!\!\in {^*\mathbb{R}} \, \forall y \,^*\!\!\in {^*\mathbb{R}} \ \ {^*f}(x \,^*\!\!+ y) = {^*f}(x) \,^*\!\!+ {^*f}(y).$$
We often omit * from $^*f$ and $^*+$.

\end{enumerate}
\end{example}

\begin{definition}(Transfer Principle II: *-Property, *-Finitenss, *-Finite Sum)
\begin{enumerate}
\item Let $P (\in U)$ define some property $Prop$. We say $u$ is $Prop$ if $u \in P$. In this situation, we say $v$ is *-$Prop$ if $v \,^*\!\in \,^*P$. An example is the following.  For $A \in U$, if $P \equiv \mathcal{P}_F(A)$ (the set of all the finite subsets of $A$), then $u (\in P)$ is a finite subset of $A$ and $v (\,^*\!\in \,^*P)$ is a *-finite subset of $^*A$.
\item Denote the finite sum of the elements of a finite subset of $\mathbb{R}$ by $\Sigma: \mathcal{P}_F(\mathbb{R}) \mapsto \mathbb{R}$. Then, we obtain the *-finite sum $^*\Sigma: \,^*\mathcal{P}_F(\mathbb{R}) \mapsto \ ^*\mathbb{R}$. That is, *-finite sum is defined on all the *-finite subset of $^*\mathbb{R}$.
\end{enumerate}
\end{definition}
 
\begin{example}(Concurrence Principle)
\begin{enumerate}
\item Since the formula $\phi(a,b) \equiv a<b \land b \in \mathbb{R}$ is concurrent with respect to $\mathbb{R}$, we obtain $b \,^*\!\in \,^*\mathbb{R}$ such that for any $a \in \mathbb{R}$ $^*a {\,^*\!\!<}  b$. That is, $b$ is an infinite number and $1/b$ is an infinitesimal. For $x, y \,^*\!\in \,^*\mathbb{R}$, we write $x \simeq y$ if $x-y$ is infinitesimal.
\item The formula $\phi(a,b) \equiv a \in b \land b \in \mathcal{P}_F([0,1])$ is concurrent with respect to $[0,1]$. Hence, there exists $b \,^*\!\in \,^*\mathcal{P}_F([0,1])$ such that $^*a \,^*\!\in b$ for all $a \in [0,1]$. In other words, there exists
*-finite subset of $^*[0,1]$ that contains all the elements of $[0,1]$, if we assume $[0,1] \subseteq \,^*[0,1]$.
\end{enumerate}
\end{example}

\begin{example}(Standard Part)
If $c {\,^*\!\in} {\,^*\mathbb{R}}$ is finite, $\sup(\{ x \in \mathbb{R} \,|\, ^*x  \,^*\!\!<  c \} )$ is called the standard part of $c$ and denoted by ${\rm st}(b)$. It is easy to check $c \simeq {\rm st}(b)$.
\end{example}

\begin{lemma}(Uniform Continuity)\label{uniformcontinuity} 
Let $f(x)$ be a function on $[0,1]$. Then  $f(x)$ is uniformly contiuous on $[0,1]$ if only if $\forall x,y {\,^*\!\in} {^*[0,1]} \  (x \simeq y \to f(x) \simeq f(y))$.
\end{lemma}
\begin{proof}
Suppose that $f(x)$ is uniformly contiuous on $[0,1]$. Then by definition $\forall \epsilon \in \mathbb{R}_+ \exists \delta \in \mathbb{R}_+ \ (|x-y|<\delta \to |f(x)-f(y)|<\epsilon)$. Fix any $\epsilon \in \mathbb{R}_+$. For $x,y \,^*\!\in {^*[0,1]}$, if $x \simeq y$ then obviously $|x-y| <\delta$. By the Transfer Principle $|f(x)-f(y)|<\epsilon$ so that $f(x)\simeq f(y)$ because $\epsilon \in \mathbb{R}_+$ is arbitary. Conversely, suppose that $\forall x,y {\,^*\!\in} {^*[0,1]}\  (x \simeq y \to f(x) \simeq f(y))$. Fix any $\epsilon \in \mathbb{R}_+$. If $|x-y|$ is less than some positive infinitesimal, $|f(x)-f(y)|<\epsilon$ by assumption. That is $\exists \delta {\,^*\!\in} {\,^*\mathbb{R}_+} (|x-y|<\delta \to |f(x)-f(y)|<\epsilon)$. By the Transfer Principle we obtain $\exists \delta \in \mathbb{R}_+ (|x-y|<\delta \to |f(x)-f(y)|<\epsilon)$.
\end{proof}

\begin{definition}(Good *-Partition of $[0,1]$) 
$p \equiv \{0=a_0<a_1<\cdots<a_i<\cdots<a_n =1\} \ (a_i \in \mathbb{R},n \in \mathbb{N})$ is a partition of $[0,1]$. By the Transfer Principle, $P\equiv \{0=a_0<a_1<\cdots<a_i<\cdots<a_N =1\} \ (a_i \in {^*\mathbb{R}}, N \in {^*\mathbb{N}})$
 is a *-partition of $^*[0,1]$. $P$ is called "good" if $\mathbb{R} \subseteq P$ The term "good" is used only in the next example. 
\end{definition}

\begin{example}(Riemann-Stieltjes Integral)
Let $f(x): [0,1] \mapsto \mathbb{R}$ be a continuous function and $g(x): [0,1] \mapsto \mathbb{R}$ be a non-decreasing function. For $p \equiv \{0=a_0<a_1<\cdots<a_n =1\}$ (a partition of $[0,1]$), set $ S(p) \equiv \sum_{k=1}^n f(a_{k-1})(g(a_k)-g(a_{k-1}))$. $S$ is a function from the set of all the partition of [0,1] $I$ to $\mathbb{R}$ so that by the Transfer Principle, $^*S: {^*I} \mapsto {^*\mathbb{R}}$. In other words, for a *-partition of $^*[0,1], P \equiv \{0=a_0 < a_1< \cdots <a_N=1\}$ $^*S(P) = {^*\sum}_{k=1}^N f(a_{k-1})(g(a_k) -g(a_{k-1}))$, where $N \in {^*\mathbb{N}}$ and $^*\sum$ is *-finite sum. Note that *'s are omitted from $^*f, ^*g$ and $^*-$. Suppose that $P$ and $P'$ are "good" *-partitions of $^*[0,1]$. Then $^*S(P) \simeq {^*S(P')}$. To see this, let $P''$ be the combined *-partitions of $P$ and $P'$. Then, it is a routine to verify that $S(P) \simeq S(P'')$ and $S(P') \simeq S(P'')$ by using Lemma \ref{uniformcontinuity}.  
\end{example}


\begin{thebibliography}{9}
\bibitem{Bernstein}
Bernstein AR, 
The Spectral Theorem - a Non-Standard Approach,     
Mathematical Logic Quarterly, {\bf18} (1972), 419-434.
\bibitem{Davis}
Davis M,   
Applied nonstandard analysis,
John Wiley \& Sons (1977), Dover Publications (2005).
\bibitem{Goldbring}
Goldbring I, 
A nonstandard proof of the spectral theorem for unbounded self-adjoint operators, 
Expositiones Mathematicae, {\bf39} (2021), 590-603.
\bibitem{Halmos}
Halmos PR, 
Finite-dimensional vector spaces, 2nd edition, 
Van Norstrand (1958), Dover Publications (2017).
\bibitem{Moore}
Moore LC Jr.,
Hyperfinite extensions of bounded operators on a separable Hilbert space,     
Transactions of the American Mathematical Society, {\bf218} (1976), 285-295.
\bibitem{Raab}
Raab A, 
An approach to nonstandard quantum mechanics,
Journal of Mathematical Physics, {\bf 45}(2004), 4791-4809.
\bibitem{Yamashita}
Yamashita H and Ozawa M, 
Nonstandard representations of unbounded self-adjoint operators,
RIMS Kokyuroku, {\bf1186}(2001), 106-118, https://www.kurims.kyoto-u.ac.jp /\~{}kyodo/kokyuroku/contents /pdf/1186-11.pdf
\end{thebibliography}
\end{document}